\documentclass[12pt,a4paper]{article}
\usepackage[english]{babel}
\usepackage{amsmath}
\usepackage{amstext}
\usepackage{amssymb}
\usepackage{amsthm}
\usepackage{amsfonts}
\usepackage{fullpage}
\usepackage{graphicx}
\usepackage{epsf}
\usepackage[dvips]{epsfig}
\usepackage{stmaryrd}
\usepackage[T1]{fontenc}
\usepackage{latexsym}
\usepackage{psfrag}
\usepackage{cyrillic}
\usepackage{hyperref}
      
\title{Finite blocking property {\em vs} pure periodicity}

\author{Thierry Monteil \footnote{CNRS -- 
LIRMM,
161 rue Ada,
34392 Montpellier,
France -- 
\url{http://www.lirmm.fr/~monteil/}
}}

\date{}

\newcommand{\N}{\mathbb N}
\newcommand{\R}{\mathbb R}
\newcommand{\C}{\mathbb C}
\newcommand{\Z}{\mathbb Z}
\newcommand{\Q}{\mathbb Q}
\newcommand{\defi}{\stackrel{\mbox {\tiny def}}{=}}

\newenvironment{appli}{\left( \begin{array}{ccc}}{\end{array} \right)}
\newcommand{\dans}{& \longrightarrow &}
\newcommand{\donne}{& \longmapsto &}
\newcommand{\ba}{\begin{appli}}
\newcommand{\ea}{\end{appli}}
\newtheorem{lem}{Lemma}[]
\newtheorem{prop}{Proposition}[]

\newenvironment{demo}{\noindent {\bf Proof:}}{\hfill $\square$ \\}
\newtheorem{thm}{Theorem}[]

\newcommand{\B}{\mathbb B} 
\renewcommand{\H}{\mathcal H}
\newcommand{\resp}{resp.\ }

\newcommand{\T}{{\mathcal{T}}}

\newcommand{\card}{\mbox{card}}

\renewcommand{\P}{{\mathcal{P}}}
\renewcommand{\S}{{\mathcal{S}}}

\renewcommand{\C}{{\mathcal{C}}}
\renewcommand{\B}{{\mathcal{B}}}

\renewcommand{\T}{{\mathcal{T}}}

\begin{document}

\maketitle


\begin{abstract}
\noindent A translation surface $\S$ is said to have the finite
blocking property if for every pair $(O,A)$ of points in $\S$
there exists a finite number of ``blocking'' points $B_1, \dots ,
B_n$ such that every geodesic from $O$ to $A$ meets one of the
$B_i$'s. $\S$ is said to be purely periodic if the directional
flow is periodic in each direction whose directional flow contains
a periodic trajectory (this implies that $\S$ admits a cylinder
decomposition in such directions). We will prove that finite
blocking property implies pure periodicity. We will also classify
the surfaces that have the finite blocking property in genus $2$:
such surfaces are exactly the torus branched coverings. 
Moreover, we prove that in every stratum, such surfaces form a
set of null measure. In the Appendix, we prove that completely periodic translation surfaces form a set of null measure in every stratum.\\



\end{abstract}

\clearpage

\section*{Introduction}

A {\em translation surface} is a triple $(\S, \Sigma, \omega)$
such that $\S$ is a topological compact connected surface,
$\Sigma$ is a finite subset of $\S$ (whose elements are called
{\em singularities}) and $\omega = (U_i,\phi_i)_{i\in I}$ is an
atlas of $\S \setminus \Sigma$ (consistent with the topological
structure on $\S$) such that the transition maps (i.e. the $\phi_j
\circ \phi_i^{-1} : \phi_i(U_i\cap U_j) \rightarrow \phi_j(U_i\cap
U_j)$ for $(i,j)\in I^2$) are translations. This atlas gives to
$\S \setminus \Sigma$ a Riemannian structure; therefore, we have
notions of length, angle, measure, geodesic... Moreover, we assume
that $\S$ is the completion of $\S \setminus \Sigma$ for this
metric.

A translation surface can also be seen as a holomorphic
differential $h$ on a Riemann surface (the singularities
correspond to the zeroes of the differential, and in an admissible
atlas $\omega$, $h$ is of the form $h=dz$). This point of view is
useful to give coordinates (and therefore a topology and a
measure) to the moduli space of all translation surfaces,
while the first one allows to make pictures in the plane.

Translation surfaces provide one of the main tool for the study of rational polygonal billiards.\\ \\

Since the unit tangent bundle of $\S$ enjoys a canonical global
decomposition $U\S=\S \times \mathbb{S}^1$
, the study of the geodesic flow on $\S$ can be done
through two points of view depending on whether the variable is
the first or the second projection:

\begin{enumerate}

\item We can fix one particular direction $\theta\in\mathbb{S}^1$:
this corresponds to the study of the {\em directional flow}
$\phi_{\theta} : \S\times \R \rightarrow \S $. We are usually
interested in the ergodic properties of this flow depending on
whether $\theta$ is a saddle connection direction or not, or on
the way $\theta$ can be approximated by such saddle connections
(see \cite{KMS}, \cite{Ve}, \cite{MT}, \cite{Vo}, \cite{Ma}, \cite{MS}, \cite{Ch}).
This class of problems is called {\em dynamical problems}.
This point of view is the most studied by mathematicians.

\item We can also fix one point in $x\in\S$: this corresponds to
the study of the {\em exponential flow} $exp_x : \mathbb{S}^1
\times \R \rightarrow \S $. We are usually interested in looking
at which points we can reach when we let the second variable vary,
or at the statistics of the return to (or near) $x$ (see \cite{Kl},
\cite{To}, \cite{Mo-HS}, \cite{Mo-V}, \cite{BGZ}, \cite{CG}).
This class of problems is called {\em illumination problems}.
This point of view is connected with computer science
(computational geometry, ray tracing, art gallery).

\end{enumerate}

Note that those two flows are defined on a dense $G_\delta$-set of
full measure but not everywhere because the singularities don't
allow some (few) geodesics to be defined on $\R$.\\


We are interested here by the connection between an illumination
problem called the finite blocking property, and a dynamical one
named the pure periodicity.\\

A translation surface $\S$ is said to have the {\em finite
blocking property} if for every pair  $(O,A)$ of points in $\S$,
there exists a finite number of points $B_1, \dots , B_n$
(different from $O$ and $A$) such that every geodesic from $O$ to
$A$ meets one of the $B_i$'s.\\

A translation surface $\S$ is said to be {\em purely periodic} if
for any $\theta\in\mathbb{S}^1$, the existence of a (non-singular) periodic orbit
in the direction $\theta$ implies that the directional flow
$\phi_{\theta}$ is periodic (i.e. there exists $T>0$ such that
$\phi_{\theta}^T=Id_{\S}$ a.e.). 
In particular, the widths of the cylinders in such a direction $\theta$ are commensurable. 
This notion is a stronger version of the complete periodicity defined in \cite{Cal}.\\

The paper is organized as follows:

In the first section, we recall some background and known results
about the finite blocking property. The second one is devoted to
the proof of the following result:\\

{\bf Theorem \ref{th-periodicity}.} {\em Let $\S$ be a translation
surface with the finite blocking property. Then $\S$ is purely periodic.}\\

The third section is devoted to some applications of this theorem.
First, we give an improvement of \cite{Mo-V}, Theorem 3:\\

{\bf Theorem \ref{th-measure}.} {\em In genus $g\geq 2$, the set
of translation surfaces that fail the finite blocking property is
of full measure in every stratum.}\\

Then, we give the complete classification in genus $2$:\\

{\bf Theorem \ref{th-genus2}.} {\em Let $\S$ be a translation
surface of genus $2$. Then $\S$ has the finite blocking property
if and only if $\S$ is a torus branched covering.}\\

In the Appendix, we prove the following known result which is needed in the proof of Theorem \ref{th-measure}:\\


{\bf Theorem \ref{th-mesure-zero-completement-periodiques}.}
{\em In genus $g\geq 2$, the set of completely periodic translation surfaces has measure zero in every stratum.}\\

\clearpage

\section{Background}\label{section-background} 
In this section, we present some results about the finite blocking property
that were proved in \cite{Mo-V}.

\subsection{Stability}
The finite blocking property is stable under some classical
operations on translation surfaces:

\begin{description}

\item[Singularities]

We have used the convention that a geodesic stops when it reaches
a singularity. This convention is not restrictive for the study of
the finite blocking property since we can add those singularities
in any blocking configuration. Moreover the finite blocking
property is not perturbed when we add or remove removable
singularities (a singularity $\sigma \in \Sigma$ is said to be
{\em removable} if there exists an atlas $\omega' \supset \omega$
such that $(\S, \Sigma \setminus \{ \sigma \}, \omega')$ is a
translation surface).

\item[Branched covering]

A {\em branched translation covering} between two translation surfaces is a
map $\pi : (\S,\Sigma) \rightarrow (\S',\Sigma')$ that is
a topological branched covering that locally preserves the
translation structure. A {\em torus branched covering} is a
translation surface $\S$ such that there exists a branched translation covering
from $\S$ to a flat torus.

\begin{prop}\label{prop-covering}
Let $\pi : (\S,\Sigma) \rightarrow (\S',\Sigma')$ be a translation covering branched on a finite set
$\mathcal{R}' \subset \S'$. Then $\S$ has the finite blocking property if and
only if $\S'$ has.
\end{prop}

\item[Rational billiards]

Let $\P$ denote a polygon in $\mathbb{R}^2$. Let $\Gamma \subset
O(2,\R)$ be the group generated by the linear parts of the
reflections in the sides of $\P$. When $\Gamma$ is finite, we say
that $\P$ is a {\em rational polygonal billiard table}. When $\P$ is
simply connected, $\P$ is rational if and only if all the angles
between edges are rational multiples of $\pi$. A classical
construction due to Zemljakov and Katok (see \cite{ZK}, \cite{MT})
allows us to associate to each rational billiard table $\P$ a
translation surface $ZK(\P)$
(as an example, the translation surface associated to the square
billiard table is the torus $\mathbb{T}^2 = \R^2 / \Z^2 $).

\begin{prop}\label{prop-ZK}
Let $\P$ be a rational polygonal billiard table. Then $ZK(\P)$ has the
finite blocking property if and only if $\P$ has.
\end{prop}

\item[Action of \mathversion{bold}$GL(2,\R)$ \mathversion{normal}]

If $A\in GL(2,\R)$, we can define the translation surface $$A .
(\S, \Sigma,(U_i,\phi_i)_{i\in I}) \defi (\S, \Sigma, (U_i,A \circ
\phi_i)_{i\in I});$$ hence we have an action of $GL(2,\R)$ on the
class of translation surfaces. We classically consider only
elements of $SL(2,\R)$ (see \cite{MT}), but area is not assumed to be preserved here.

\begin{prop}\label{prop-GL2R}
Let $\S$ be a translation surface and $A$ be in $GL(2,\R)$. Then
$\S$ has the finite blocking property if and only if $A.\S$ has.
\end{prop}

\end{description}

\subsection{Sufficient conditions}

\begin{prop}[Fomin \cite{Fomin}]\label{prop-Fomin}
The square billiard table has the finite blocking property.
\end{prop}

Hence, combining propositions \ref{prop-covering}, \ref{prop-ZK}, \ref{prop-GL2R}
and \ref{prop-Fomin}, we have:

\begin{prop} \label{prop-torus-branched-covering}
Any torus branched covering has the finite blocking property.
\end{prop}

\subsection{Necessary conditions}
We first give some definitions and conventions concerning isometries and cylinders.
In this paper, ``isometry'' means ``orientation-preserving local isometry''.
Since an isometry $i$ from an open subset $U$ of $\R^2$ to a translation surface $\S$ can be uniquely extended by continuity to
$\overline{U}$ (it is uniformly continuous), we will also denote by $i$ this extension (that is also an isometry).
In this process, we may lose the injectivity of $i$.\\

A {\em subcylinder} $\C$ is an isometric copy of $\R / w \Z \times
(0,h)$ in a translation surface $\S$ ($w>0$, $h>0$). $w$ and $h$
are unique and called the {\em width} and the {\em height} of $\C$.
The {\em direction} of $\C$ is the direction of the image of $\R / w \Z \times \{ h /2 \}$
(which is considered here as an oriented closed geodesic); it is defined modulo $2\pi$.
The images in $\S$ of $\R / w \Z \times \{0\}$ and  $\R / w \Z \times \{h\}$
are called the {\em sides} of $\C$.
The side $\R / w \Z \times \{h\}$ is called the {\em upper side}.

A {\em cylinder} is a maximal subcylinder (for the inclusion).
Any (regular) periodic trajectory for $\phi_\theta$ can be
thickened to obtain a cylinder in the direction $\theta$.
Whenever $\S$ is not a torus, each side of a cylinder contains at least one
singularity and is a finite union of saddle connections
(a {\em saddle connection} is a geodesic joining two singularties).
We say that $\S$ admits a {\em cylinder decomposition} in the direction $\theta$ if
the (necessarily finite) union of cylinders whose direction is $\theta$ is dense in $\S$
(then, the complement of this union is the union of the sides of those cylinders).\\

\begin{figure}[!h]
\begin{center}
\psfrag{w1}{$w_1$} \psfrag{h1}{$h_1$} \psfrag{w2}{$w_2$}
\psfrag{h2}{$h_2$} \psfrag{C1}{$\C_1$} \psfrag{C2}{$\C_2$}
\psfrag{O}{$(0,0)$} \psfrag{l}{$l$}
\includegraphics[scale = 0.56]{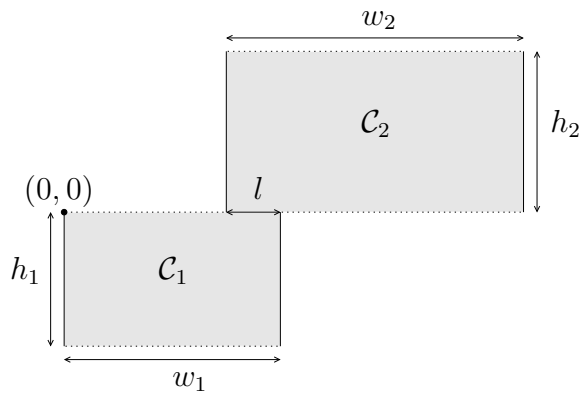}
\caption{\label{2cylinder.eps} The opposite vertical sides are
identified by translation while the dotted horizontal ones are
glued with the rest of the surface; $l>0$.}
\end{center}
\end{figure}

Let $\C_1$ and $\C_2$ be two cylinders with width $w_1$ (\resp
$w_2$) and height $h_1$ (resp. $h_2$) in a translation surface
$\S$. We suppose that $\C_1$ and $\C_2$ are two different
cylinders in the same direction whose closures have a nontrivial
(i.e. not reduced to a finite set) intersection. The situation is described in Figure \ref{2cylinder.eps}.

\begin{lem}\ \label{lem2cyln}
Let $\S$ be a translation surface that contains two different
cylinders $\C_1$ and $\C_2$ in the same direction whose closures
have a nontrivial intersection. If their widths are
uncommensurable, then $\S$ fails the finite blocking property.
\end{lem}

This lemma implies all the remaining results of section
\ref{section-background} and will be useful in the proof of
Theorem \ref{th-periodicity}.

\begin{prop}
The only regular polygons (considered as billiard tables) with the
finite blocking property are the equilateral triangle, the square
and the regular hexagon.
\end{prop}

\begin{prop}\label{prop-Veech}
A Veech surface has the finite blocking property if and only if it
is a torus covering, branched over only one point.
\end{prop}

In other terms, arithmetic surfaces are exactly the Veech surfaces
with the finite blocking property.
Note that proposition \ref{prop-Veech} was proven
independently 
by Eugene Gutkin, and is the main result of
\cite{Gutkin}.


\subsection{Global point of view: moduli spaces of holomorphic forms}\label{subsection-moduli-space}

A singularity $\sigma \in \Sigma$ has a cone angle of the form
$2(k+1)\pi$, with $k\geq 0$; we say that $\sigma$ is of {\em
multiplicity} $k$. Removable singularities are the singularities of multiplicity $0$.
In terms of holomorphic differentials, it is equivalent to say that
there is a chart around $\sigma$ such that $h=z^{k}dz$ (i.e. $\sigma$ is a zero of order $k$ of $h$).

If $1\leq k_1\leq k_2 \leq \dots \leq k_n$ is
a sequence of integers whose sum is even, we denote by $\H(k_1, k_2,
\dots, k_n)$ the {\em stratum} of translation surfaces with
exactly $n$ singularities whose multiplicities are
$k_1, k_2, \dots, k_n$ (we consider only surfaces without removable singularities). A translation surface in $\H(k_1, k_2,
\dots, k_n)$ has genus $g=1+(k_1+k_2+ \dots +k_n)/2$.

Each stratum carries a natural topology and a $SL(2,\R)$-invariant measure that are 
defined in \cite{Masur-ei} and \cite{Veech-ei}. Let $\S$ be an element of some $\H(k_1, k_2, \dots, k_n)$ and 
let $\B=\{\gamma_1, \dots , \gamma_{2g+n-1}\}$ be a basis of the {\em relative} homology of $\S$:
it is just the concatenation of a basis of the first topological homology group $H_1(\S,\Z)$ 
with a set of $n-1$ curves from a singularity to the other ones.
If $\S'$ is another translation surface (built on the same topological surface), 
let us denote by $hol_{\S'}$ the associated (translational) holonomy.
The map 
$$\Phi \defi \ba \H(k_1, k_2,\dots, k_n)  \dans \mathbb{C}^{2g+n-1} \simeq (\R^2)^{2g+n-1} \\ 
\S'  \donne (hol_{\S'}(\gamma_1) , \dots , hol_{\S'}(\gamma_{2g+n-1}))  \ea$$
is named the {\em period map} and is a local homeomorphism in a neighbourhood of $\S$ 
and in this system of coordinates, the former measure is absolutely continuous relatively to Lebesgue.\\




\begin{prop}\label{prop-dense}
In genus $g\geq 2$, the set of translation surfaces that fail the
finite blocking property is dense in every stratum.
\end{prop}

In fact, we proved in \cite{Mo-V} that the hypotheses of Lemma \ref{lem2cyln} are satisfied on a dense subset of each stratum.
Unfortunately, those hypotheses are valid only on a set of zero measure.
Proposition \ref{prop-dense} will be considerably improved by Theorem \ref{th-measure}.

\clearpage

\section{Main result}

A translation surface is said to be {\em completely periodic} if the
following property holds (see \cite{Cal}): if there exists a periodic
orbit in the direction $\theta$, then the surface admits a cylinder
decomposition in the direction $\theta$.

\begin{prop}\label{prop-complete-periodicity}
Let $\S$ be a translation surface with the finite blocking
property. Then $\S$ is completely periodic.
\end{prop}

\begin{demo} Let $\theta$ be a direction that contains a periodic
orbit. Up to a rotation, we assume
that this direction is the horizontal one (defined by $\theta=0 \mod 2\pi$). We can thicken all of
such horizontal periodic orbits to obtain a finite number of cylinders. Assume by
contradiction, that $\S$ is not fully decomposed into cylinders in
the horizontal direction: there exists a cylinder $\C$ and singularities $\sigma$ and $\sigma'$ (possibly equal) that lie on the upper side of $\C$
such that the horizontal geodesic $\gamma$ (whose length is denoted by $l(\gamma)$ and image by $|\gamma|$)
from $\sigma$ to $\sigma'$ exists (the singularities are consecutive for the cyclic order on the upper side of $\C$)
and is not a part of a side of any other cylinder.
To be more precise, there exists $s>0$ (by compactness of $|\gamma|$ and by finiteness of the number of horizontal cylinders)
such that there is an injective isometry $i: (0,l(\gamma)) \times  (0,s) \rightarrow \S$ such that $i([0,l(\gamma)]\times\{0\})=|\gamma|$
and $i( (0,l(\gamma))\times (0,s))$ doesn't meet any horizontal cylinder.
Up to a dilatation (that belongs to $GL^{+}(2,\R)$), we assume that the width of $\C$ is $1$ and that its height is $2$.
The situation is described in Figure \ref{cyln-minim.eps}.

\begin{figure}[!h]
\begin{center}
\psfrag{g}{$|\gamma|$}
\psfrag{1}{$1$}
\psfrag{2}{$2$}
\psfrag{lg}{$l(\gamma)$}
\psfrag{C}{$\C$}
\psfrag{s}{$s$}
\psfrag{S}{$\sigma$}
\psfrag{S'}{$\sigma'$}
\includegraphics[scale = 0.58]{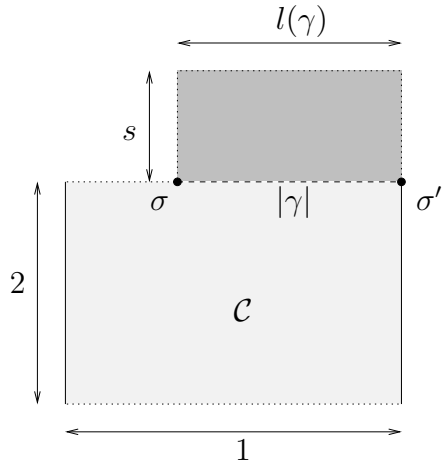}
\caption{\label{cyln-minim.eps} The opposite vertical black sides are
identified by translation; the dotted ones are glued with the rest of the surface;
the dark zone doesn't meet any horizontal cylinder.}
\end{center}
\end{figure}

The heart of the proof is organised into three steps. A summary
concludes each of them.

\begin{description}

\item[Step 1] {\em Construction of certain thin rectangles in $\S$ that contain a common point $A$.}

Let $\varepsilon > 0$ be smaller than $s$.

Let $T_{\varepsilon}$ be the supremum of all $t>0$ such that there exists an injective isometry
$j_t : (0,t) \times (0,\varepsilon) \rightarrow \S$ that coincides with $\gamma$ on $(0,\min(l(\gamma),t)) \times\{0\}$.
$T_{\varepsilon}$ is greater than or equal to $l(\gamma)$ and therefore strictly positive.
$\S$ has finite area hence $T_{\varepsilon}$ is finite.
Moreover $T_{\varepsilon}$ is a maximum (indeed, if $0 < t \leq t' < T_{\varepsilon}$, then $j_{t'}$ extends $j_t$):
there exists an injective isometry $j_{T_\varepsilon} : (0,T_{\varepsilon}) \times (0,\varepsilon) \rightarrow \S$
that coincide with $\gamma$ on $(0,l(\gamma)) \times\{0\}$.

There are only three possibilities depending on what happens on $j_{T_\varepsilon}(\{T_{\varepsilon}\}\times (0,\varepsilon))$
(if $j_{T_\varepsilon}(\{T_{\varepsilon}\}\times(0,\varepsilon))$ does not meet a singularity,
then $j_{T_\varepsilon}$ can be extended to an isometry that fails to be injective,
hence $I_{\varepsilon}=j_{T_\varepsilon}(\{T_{\varepsilon}\}\times (0,\varepsilon))$
and $J_{\varepsilon}=j_{T_\varepsilon}(\{0\}\times (0,\varepsilon))$ intersect each other):

\begin{enumerate}

\item $I_{\varepsilon}$ meets a singularity denoted by $\sigma_{\varepsilon}$.

\item $I_{\varepsilon} = J_{\varepsilon} $:
we just have constructed a cylinder, which is forbidden by the assumption, so this case does not occur.

\item The vertical intervals $I_{\varepsilon}$ and $J_{\varepsilon}$
intersect each other but we are neither in the first nor in the second case. Since $\sigma=j_{T_\varepsilon}(0,0)$ is a singularity,
$I_\varepsilon \cap J_{\varepsilon}$ is a subinterval of $J_{\varepsilon}$ that hits the extremity of $J_{\varepsilon}$ that is not $\sigma$
(i.e. $j_{T_\varepsilon}(0,\varepsilon)$) (see Figure \ref{enroulement.eps}).


\begin{figure}[!h]
\begin{center}
\psfrag{g}{$|\gamma|$}
\psfrag{1}{$1$}
\psfrag{2}{$2$}
\psfrag{e}{$\varepsilon$}
\psfrag{C}{$\C$}
\psfrag{s}{$s$}
\psfrag{S}{$\sigma$}
\psfrag{S'}{$\sigma'$}
\psfrag{O}{$O$}
\psfrag{A}{$A$}
\psfrag{h}{$h$}
\psfrag{l}{$l$}
\psfrag{e'}{$\varepsilon'$}
\psfrag{Ie}{$I_{\varepsilon}$}
\psfrag{Je}{$J_{\varepsilon}$}
\includegraphics[scale = 0.58]{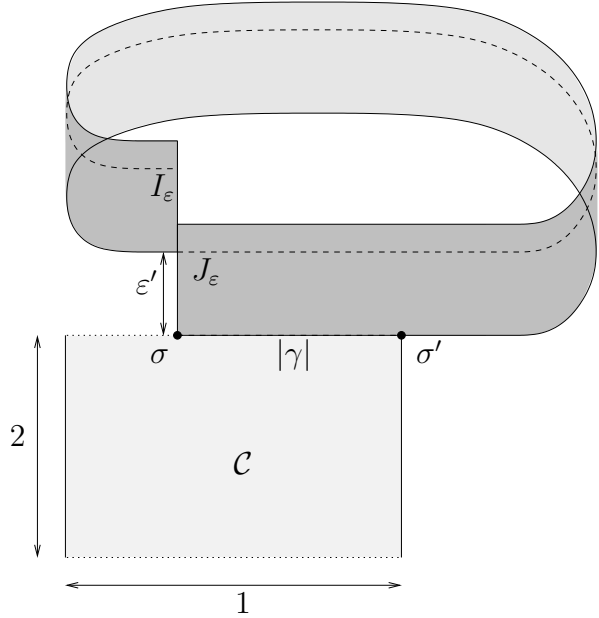}
\caption{\label{enroulement.eps} The dark zone is the image of $j_{T_\varepsilon}$.}
\end{center}
\end{figure}

So, there exists $0 < \varepsilon' \leq \varepsilon$ such that $J_{\varepsilon} \setminus I_{\varepsilon} =  j_{T_\varepsilon}(\{0\}\times (0,\varepsilon'])$.
Now, the same three possibilities are offered to $j_{T_{\varepsilon'}}$, but the third case can't happen. Indeed, otherwise there is $0<\varepsilon'' < \varepsilon'$ such that $j_{T_{\varepsilon'}}(0,\varepsilon')=j_{T_{\varepsilon'}}(T_{\varepsilon},0)=j_{T_{\varepsilon'}}(T_{\varepsilon'},\varepsilon'-\varepsilon'')$ is not a singularity but 
coming out of it are at least 2 horizontal lines
(equivalently, this point has a conical angle greater than $2\pi$).


\end{enumerate}

So, up to reducing the size of $\varepsilon$ (by replacing $\varepsilon$ by $\varepsilon'$ whenever we are in the third case),
we can consider that we are in the first case, and,
since $\Sigma$ is finite, there exists a sequence $(\varepsilon_n)_{n\in\N}$ of strictly positive numbers that converges to $0$
and a singularity $\sigma_{\infty} \in \Sigma$ such that for all $n$ in $\N$, $\sigma_{\varepsilon_n}=\sigma_{\infty}$.
If $n\in \N$, there exists $0<y_n<\varepsilon_{n}$ such that $\sigma_{\infty}=j_{T_{\varepsilon_{n}}}(T_{\varepsilon_n}, y_n)$.
Let $A=j_{T_{\varepsilon_{n}}}(T_{\varepsilon_n} - l(\gamma)/2 , y_n)$
($A$ does not depend on $n\in\N$ and is in the image of every $j_{T_{\varepsilon_n}}$). 

{\em
To sum up this step and fix notations, we proved that there exists
$A\in \S$ such that for all $\varepsilon > 0$ there exists
$l_{\varepsilon} \geq l(\gamma)$ and
$0 < h_{\varepsilon} \leq \varepsilon$
such that there exists an injective isometry $i_{\varepsilon}$ from
$(0,l_{\varepsilon}) \times (0,h_{\varepsilon})$ to $\S$
that coincides with $\gamma$ on $(0,l(\gamma)) \times\{0\}$
and whose image contains $A$.}

The situation is summarized in Figure \ref{situation-step-1.eps}. Let $O$ be the point drawn on this Figure.
We denote $A=i_{\varepsilon}(x_{\varepsilon}, y_{\varepsilon})$.

\begin{figure}[!h]
\begin{center}
\psfrag{g}{$|\gamma|$}
\psfrag{1}{$1$}
\psfrag{e}{$\varepsilon$}
\psfrag{C}{$\C$}
\psfrag{s}{$s$}
\psfrag{S}{$\sigma$}
\psfrag{S'}{$\sigma'$}
\psfrag{O}{$O$}
\psfrag{A}{$A$}
\psfrag{h}{$h_{\varepsilon}$}
\psfrag{l}{$l_{\varepsilon}$}
\psfrag{y}{$y_{\varepsilon}$}
\includegraphics[scale = 0.6]{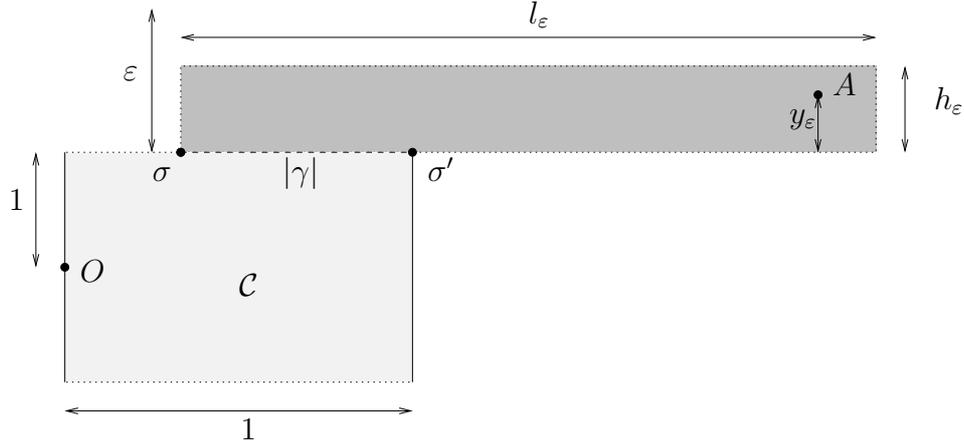}
\caption{\label{situation-step-1.eps} 
Summary of Step 1. }
\end{center}
\end{figure}

\item[Step 2] {\em Construction of ``teepees'' of more and more geodesics between $O$ and $A$.}

Let $\varepsilon>0$.
The unfolding of $\C$ allows us to define an extension of $i_{\varepsilon}$ as in Figure \ref{thales1.eps}:
we get an isometry $i_{\varepsilon} : (0,l_{\varepsilon}) \times [0,h_{\varepsilon}) \bigcup \R\times (-2,0) \rightarrow \S$.
In $\R^2$, we denote without ambiguity $\sigma=(0,0)$, $\sigma'=(l(\gamma),0)$ and $A=(x_{\varepsilon}, y_{\varepsilon})$.

\begin{figure}[!h]
\begin{center}
\psfrag{g}{$|\gamma|$}
\psfrag{1}{$1$}
\psfrag{e}{$\varepsilon$}
\psfrag{C}{$\C$}
\psfrag{s}{$s$}
\psfrag{S}{$\sigma$}
\psfrag{S'}{$\sigma'$}
\psfrag{s}{$S_{\varepsilon}$}
\psfrag{s'}{$S'_{\varepsilon}$}
\psfrag{O}{$O$}
\psfrag{A}{$A$}
\psfrag{h}{$h_{\varepsilon}$}
\psfrag{lg}{$l(\gamma)$}
\psfrag{d}{$d(S_{\varepsilon},S_{\varepsilon}')$}
\psfrag{y}{$y_{\varepsilon}$}
\includegraphics[scale = 0.44]{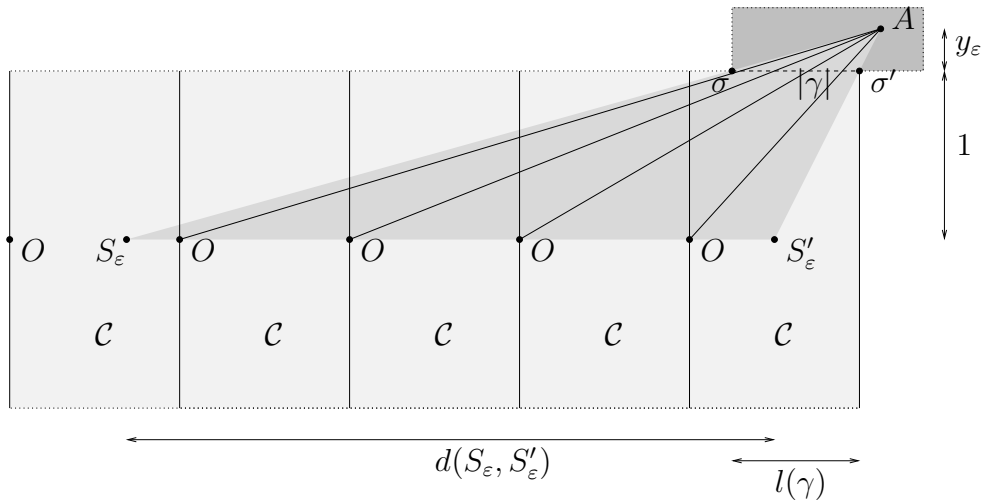}
\caption{\label{thales1.eps} First (topological) application of the theorem of Thales.}
\end{center}
\end{figure}

Let $S_{\varepsilon}$ (resp. $S_{\varepsilon}'$) be the point of the line that passes through $A$
and $\sigma$ (resp. $\sigma'$) whose y-coordinate is $-1$.
Since $(\sigma, \sigma')$ and $(S_{\varepsilon}, S_{\varepsilon}')$ are parallel,
the theorem of Thales (about similar triangles) asserts that $$\frac{d(S_{\varepsilon},S_{\varepsilon}')}{ (1+y_{\varepsilon})} = \frac{l(\gamma)}{ y_{\varepsilon}} , $$
where $d$ stands for the euclidean distance of $\R^2$.
Let $\mathcal{O}_{\varepsilon} = i_{\varepsilon}^{-1}(\{O\})\cap (S_{\varepsilon}, S_{\varepsilon}')$.
We have $$\card \, \mathcal{O}_{\varepsilon} \geq d(S_{\varepsilon},S_{\varepsilon}') -1 = l(\gamma) (1+y_{\varepsilon}) / y_{\varepsilon} -1  \xrightarrow[ \varepsilon \rightarrow 0]{} \infty .$$

Since the open triangle $(A,S_{\varepsilon},S_{\varepsilon}')$ is included in the domain of $i_{\varepsilon}$, we define
a set $\T_{\varepsilon}$ (called a {\em teepee}) of geodesics in $\S$ from $O$ to $A$ by taking the images by $i_{\varepsilon}$
of the segments joining an element of $\mathcal{O}_{\varepsilon}$ to $A$.
Of course, $$\card \, \T_{\varepsilon} = \card \, \mathcal{O}_{\varepsilon}  \xrightarrow[ \varepsilon \rightarrow 0]{} \infty  . $$

{\em
To sum up this step, we constructed, for all $\varepsilon > 0$, a set $\T_{\varepsilon}$
of more and more geodesics from $O$ to $A$ whose images are included in
$\C \cup | \gamma | \cup i_{\varepsilon} ((0,l_{\varepsilon}) \times (0,h_{\varepsilon}))$.
}

\item[Step 3] {\em The density of the geodesics of a teepee that are blocked by a fixed point vanishes.}

Let $B$ be a fixed point in $\S$ that is different from $O$ and $A$.
For $\varepsilon >0$, let $\T_{\varepsilon} (B)$ be the set of elements of $\T_{\varepsilon}$ that are blocked by $B$
(i.e. that pass through $B$).

\begin{itemize}
\item Suppose that $B$ is not in $\C$. Let $\varepsilon >0$.
If $B$ is not in the image of $i_{\varepsilon}$, then $\T_{\varepsilon} (B)$ is empty.
Otherwise, $B$ is in $i_{\varepsilon}((0,l_{\varepsilon}) \times [0,h_{\varepsilon}))$,
then the cardinal of $\T_{\varepsilon} (B)$ is at most $1$
(indeed, the elements of $\T_{\varepsilon}$ are disjoint on $i_{\varepsilon}((0,l_{\varepsilon}) \times [0,h_{\varepsilon})) \setminus \{A\}$).

Hence, $\card \, \T_{\varepsilon}(B) / \card \, \T_{\varepsilon}  \xrightarrow[ \varepsilon \rightarrow 0]{} 0$.

\item Suppose that $B$ is in $\C$.
Let $-h$ be the y-coordinate of any preimage of $B$ by any $i_{\varepsilon}$. We can assume that $h\in (0,1)$,
otherwise $\T_{\varepsilon} (B)$ is empty for any $\varepsilon >0$.
Let $\varepsilon > 0$.
Let $\beta_1$ and $\beta_2$ be two distinct elements of $\T_{\varepsilon}(B)$:
for $i\in\{1,2\}$, there exists $O_i$ in $\mathcal{O}_{\varepsilon}$ and $B_i$ in $i_{\varepsilon}^{-1}(B)$ such that
$B_i$ lies in the segment $[O_i,A]$ whose image by $i_{\varepsilon}$ is $\beta_{i}$ (see Figure \ref{thales2.eps}).

\begin{figure}[!h]
\begin{center}
\psfrag{g}{}
\psfrag{1}{$1$}
\psfrag{e}{$\varepsilon$}
\psfrag{C}{$\C$}
\psfrag{s}{$s$}
\psfrag{S}{$\sigma$}
\psfrag{S'}{$\sigma'$}
\psfrag{s}{$S_{\varepsilon}$}
\psfrag{s'}{$S'_{\varepsilon}$}
\psfrag{O1}{$O_1$}
\psfrag{O2}{$O_2$}
\psfrag{B1}{$B_1$}
\psfrag{B2}{$B_2$}
\psfrag{A}{$A$}
\psfrag{h}{$h$}
\psfrag{pe}{$p_{\varepsilon}$}
\psfrag{qe}{$q_{\varepsilon}$}
\psfrag{y}{$y_{\varepsilon}$}
\includegraphics[scale = 0.44]{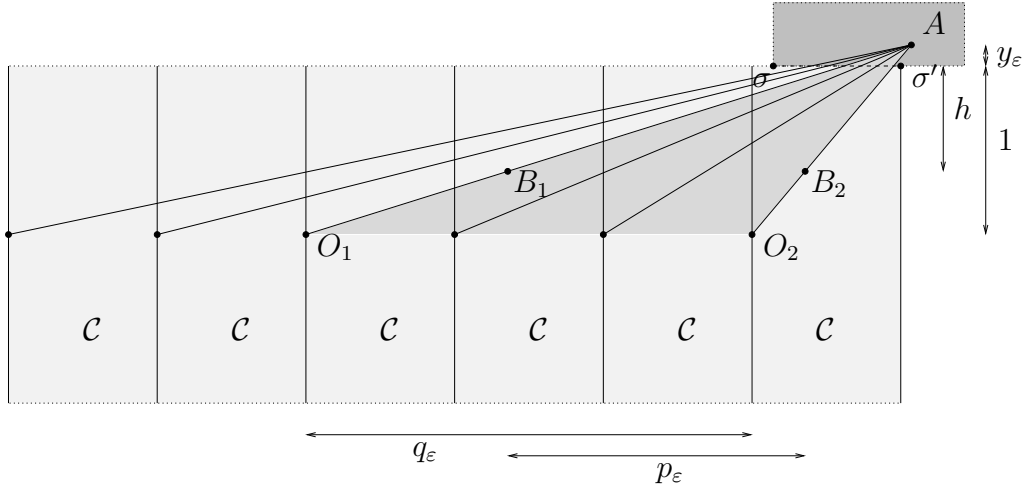}
\caption{\label{thales2.eps} Second (arithmetic) application of the theorem of Thales.}
\end{center}
\end{figure}


We choose $\beta_1$ and $\beta_2$ such that $d(O_1,O_2)$ is minimal, and set $p_{\varepsilon}=d(B_1,B_2)$ and $q_{\varepsilon}=d(O_1,O_2)$.
By minimality, $$ \card \, \T_{\varepsilon}(B) \leq \frac{\card \, \T_{\varepsilon}}{ q_{\varepsilon}} +1 . $$

Since $i_{\varepsilon}(B_1)=i_{\varepsilon}(B_2)=B$ and $i_{\varepsilon}(O_1)=i_{\varepsilon}(O_2)=O$,
$p_{\varepsilon}$ and $q_{\varepsilon}$ are 
integers.

Since $(B_1, B_2)$ and $(O_1,O_2)$ are parallel,
the theorem of Thales asserts that $$\frac{p_{\varepsilon}}{q_{\varepsilon}} = \frac{(h+y_{\varepsilon})}{ (1+y_{\varepsilon})} . $$


This quotient 
converges to $h$ when $\varepsilon$ converges to $0$ but never takes the value $h$,
hence the denominator $ q_{\varepsilon}  \xrightarrow[ \varepsilon \rightarrow 0]{} \infty$


Hence, $$\frac{\card \, \T_{\varepsilon}(B) }{ \card \, \T_{\varepsilon} } \leq \frac{1 }{ q_{\varepsilon} } + \frac{1 }{  \card \, \T_{\varepsilon}}  \xrightarrow[ \varepsilon \rightarrow 0]{} 0  . $$

\end{itemize}

{\em
To sum up this step, we proved that for any point $B$ in $\S$ that is different from $O$ and $A$,
the density of the geodesics of $\T_{\varepsilon}$ that are blocked by $B$ tends to $0$ when $\varepsilon$ tends to $0$.
}

\end{description}

Now, we assumed that $\S$ has the finite blocking property,
so there exists a finite set $\{B_1,\dots,B_n\}$ in $\S \setminus\{O,A\}$ that blocks any geodesic from $O$ to $A$.

Let $\varepsilon > 0$ be small enough such that $\Sigma_{i=1}^n \card \, \T_{\varepsilon}(B_i) / \card \, \T_{\varepsilon} < 1 $:
there exists a geodesic in $\T_{\varepsilon} \setminus \bigcup_{i=1}^{n} \T_{\varepsilon}(B_i)$ from $O$ to $A$ that doesn't meet any $B_i$,
leading to a contradiction.
So, $\S$ is completely periodic.

\end{demo}

\begin{thm}\label{th-periodicity}
Let $\S$ be a translation surface with the finite blocking
property. Then $\S$ is purely periodic.
\end{thm}

\begin{demo} Let $\theta$ be a direction that contains a periodic trajectory. By virtue of Proposition \ref{prop-complete-periodicity},
$\S$ admits a cylinder decomposition $(\C_i)_{i=1}^n$ in the direction $\theta$.
Starting from one of those cylinders and applying Lemma \ref{lem2cyln} step by step,
we show that all the widths of those cylinders are commensurable (remember that $\S$ is connected).
Hence there exists $w>0$ such that the width of $\C_i$ is equal to $p_i w$ with $p_i\in \N^*$.
Setting $T= w p_1 \dots p_n > 0$, we have $\phi_{\theta}^T = Id_{\S}$ a.e.

\end{demo}

\clearpage

\section{Some applications}

\begin{thm}\label{th-measure}
In genus $g\geq 2$, the set of translation surfaces that fail the
finite blocking property is of full measure in every stratum.
\end{thm}

\begin{demo} 
Thanks to Proposition~\ref{prop-complete-periodicity}, it suffices to prove that completely periodic translation surfaces form a set of null measure in every stratum (in genus $g\geq 2$). This fact is well known in the community and should be attributed to the folklore. 
Since it is written nowhere, we will give a proof in the Appendix.
The author would like to thank Erwan Lanneau, Samuel Leli\`evre, Anton Zorich and the referee for helpful discussions or comments about this result.

Note that Theorem \ref{th-measure} will be considerably improved in \cite{Mo-homology}.

\end{demo}

\begin{thm}\label{th-genus2}
Let $\S$ be a translation surface of genus $2$. Then $\S$ has the
finite blocking property if and only if $\S$ is a torus branched
covering.
\end{thm}

\begin{demo}

By Proposition \ref{prop-torus-branched-covering}, if $\S$ is a torus
branched covering, then $\S$ has the finite blocking property.

Conversely, suppose that $\S$ has the finite blocking property. By
Proposition \ref{prop-complete-periodicity}, $\S$ is completely
periodic. Such surfaces of genus $2$ were classified in \cite{Cal}. If
$\S$ is in $\H(2)$, then $\S$ is a Veech surface and by Proposition
\ref{prop-Veech}, $\S$ is a torus branched covering. If $\S$ is in
$\H(1,1)$, we assume by contradiction that $\S$ is not a torus branched
covering. A direct consequence of Theorem 1.2 of \cite{Cal} is that after rescaling $\S$, there
exists a square-free integer $d>0$, four positive numbers
$w_1$, $w_2$, $s_1$, $s_2$ in
$\Q(\sqrt{d})$ and a direction $\theta \in \mathbb{S}^1$ such that:



\begin{itemize}
   \item $\S$ is decomposed into cylinders $\C_1$, $\C_2$ (and sometimes $\C_3$) in the direction $\theta$
   \item $w_1$ (\resp $w_2$) is the width of $\C_1$ (\resp $\C_2$)
   \item $$\frac{w_1}{w_2}  \overline{s}_1 + \overline{s}_2 = 0$$ \\
\end{itemize}

Since $w_1$, $w_2$, $s_1$, $s_2$ are in $\Q(\sqrt{d})$, there exists two rational
numbers $r_1$ and $r_2$ such that $$\frac{w_1}{w_2} s_1 + s_2 = r_1 + r_2 \sqrt{d}.$$

The addition of those two equalities gives $\frac{w_1}{w_2} (s_1 +
\overline{s}_1) + (s_2 + \overline{s}_2) = r_1 + r_2 \sqrt{d}$.
Theorem~\ref{th-periodicity} asserts that $\S$ is purely periodic, so
$\frac{w_1}{w_2}$ is rational, so $\frac{w_1}{w_2} (s_1 +
\overline{s}_1) + (s_2 + \overline{s}_2)$ is rational, so $r_2=0$.
The subtraction of the two equalities leads to $r_1=0$.

Hence, $\frac{w_1}{w_2} s_1 + s_2$ is both null and positive,
which is impossible: $\S$ is a torus branched covering.

\end{demo}

Of course,
the notion of finite blocking property and the results presented in this paper can be translated in the vocabulary of quadratic differentials.

\clearpage

\section*{Conclusion}
We conclude this paper by the following open question: 
\begin{center} Does pure periodicity imply being a torus branched covering? \end{center}

It was positively answered in \cite{Mo-V} for Veech surfaces and in the present paper for surfaces of genus two. 
Proving that in full generality would provide an equivalence between the following assertions: 

\begin{enumerate}
              \item $\S$ is a torus branched covering
              \item there exists $A\in GL(2,\R)$ such that $hol(H_1(A.\S,\Z))=\Z + i \Z$ 
              \item $Vect_{\Q}(hol(H_1(\S,\Z)))$ has dimension $2$
              \item $\S$ has the finite blocking property
              \item $\S$ has the bounded blocking property: there exists $n\in\N$ such that for every pair $(O,A)$ of points in $\S$, there exists $n$ points $B_1, \dots , B_n$ (different from $O$ and $A$) such that every geodesic from $O$ to $A$ meets one of the $B_i$'s
              \item $\S$ has the middle blocking property: For every pair $(O,A)$ of points in $\S$, the set of the midpoints of the geodesics from $O$ to $A$ is finite
              \item $\S$ has the bounded middle blocking property: there exists $n\in\N$ such that the set of the midpoints of the geodesics from $O$ to $A$ has cardinal less than $n$
              \item $\S$ is purely periodic: in every direction $\theta$ that contains a periodic orbit, the directional flow $\phi_{\theta}$ is periodic (i.e. there exists $T>0$ such that $\phi_{\theta}^T=Id_{\S}$ a.e.)
              \item in every direction that contains a periodic orbit, $\S$ is decomposable into cylinders whose widths are commensurable
\end{enumerate}


Properties $1$, $2$ and $3$ are geometric, whereas properties $4$, $5$, $6$, $7$ are exponential and properties $8$ and $9$ are dynamic.\\

If $\C$ denotes the subgroup of $H_1(\S,\Z)$ generated by the periodic orbits, 
it is already possible to prove that if $\S$ is purely periodic, 
then  $Vect_{\Q}(hol(\C))$ has dimension $2$. 
Consequently, when $\C=H_1(\S,\Z)$, 
then the nine previous assertions are equivalent \cite{Mo-homology}.

\section{Appendix}
\begin{thm}\label{th-mesure-zero-completement-periodiques}
In genus $g\geq 2$, the set of completely periodic translation surfaces has measure zero in every stratum.
\end{thm}

\begin{demo}

Let $\S$ be a completely periodic translation surface of genus $g\geq 2$ in a stratum $\H(k_1, k_2, \dots, k_n)$. Let $\theta$ be a direction in which $\S$ admits a cylinder decomposition $(\C_i)_{i=1}^K$.

\begin{description}

\item[Step 1] {\em $\S$ can not have more than $2g+n-3$ cylinders in the direction $\theta$.}

The total conical angle around the singularities of $\S$ is equal to $\sum_{i=1}^n 2(k_i+1)\pi$.

Each side of each cylinder in the direction $\theta$ contains at least one singularity, therefore it contributes to at least $\pi$ in the total conical angle, hence any cylinder contributes to at least $2\pi$ in the total conical angle.

Now suppose that the upper side of a cylinder $\mathcal{C}_i$ meets the lower side of a cylinder $\mathcal{C}_j$ (possibly the same) along a saddle connection and suppose that each of those two sides contribute only to $\pi$ in the total conical angle. Then, each side consists in only one saddle connection from a singularity to itself, so the two sides match perfectly: this contradicts the maximality of the (sub)cylinder $\mathcal{C}_i$ (and the singularity is removable).

Hence, at least one cylinder contributes to more than $2\pi$ in the total conical angle, so the number $K$ of cylinders is at most $\sum_{i=1}^n (k_i+1) - 1 = 2g+n-3$.

Note that this bound can be strengthened by studying the separatrix diagram defined in \cite{KZ} (ask Erwan Lanneau).

\item[Step 2]  {\em The saddle connections of $\S$ in the direction $\theta$ generate a subspace of dimension at least $2$ of the relative homology $H_1(\S,\Sigma,\Z)$.}

Let $SC(\S,\theta)$ denote the set of saddle connections of $\S$ in the direction $\theta$. 
For each cylinder $\C_i$, choose a geodesic $\gamma_i$ in $\C_i$ that goes from a singularity of the lower side of $\C_i$ to a singularity of its upper side. The set $G=\{\gamma_1,\dots,\gamma_K\}\cup SC(\S,\theta)$ generates $H_1(\S,\Sigma,\Z)$ since one can follow the trajectory of any curve that join two singularities in $\S$ with a succession of elements of $G$ (localize to the cylinders crossed by the curve).

Step 1 asserts that $K \leq 2g+n-3$ and we know that $\dim(H_1(\S,\Sigma,\Z))=2g+n-1$, hence $SC(\S,\theta)$ generates a subspace of dimension at least $2$ in $H_1(\S,\Sigma,\Z)$.

\item[Step 3]  {\em $\S$ belongs to a particular set 
of measure zero in $\H(k_1, k_2, \dots, k_n)$.}

Thanks to Step 2, there exists two saddle connections $\gamma_1$ and $\gamma_2$ in $SC(\S,\theta)$ that are independent in $H_1(\S,\Sigma,\Z)$. Moreover, $\det(hol_\S(\gamma_1), hol_\S(\gamma_2))=0$ since those two vectors of $\R^2$ have the same direction $\theta$.

Now, let $\mathcal{N}(k_1, k_2, \dots, k_n)$ be the set of translation surfaces $\S'$ in $\H(k_1, k_2, \dots, k_n)$ such that there exists two independent curves $\gamma_1$ and $\gamma_2$ in $H_1(\S',\Sigma,\Z)$ such that $\det(hol_{\S'}(\gamma_1), hol_{\S'}(\gamma_2))=0$: this non-trivial relation implies that $\mathcal{N}(k_1, k_2, \dots, k_n)$ has measure zero in $\H(k_1, k_2, \dots, k_n)$ (look to the period map defined in subsection \ref{subsection-moduli-space} and note that the set of pairs $(\gamma_1,\gamma_2)\in H_1(\S',\Sigma,\Z)^2$ is countable).

We proved that $\S$ belongs to $\mathcal{N}(k_1, k_2, \dots, k_n)$.

\end{description}

Hence, the set of completely periodic translation surfaces in the stratum $\H(k_1, k_2, \dots, k_n)$ is included in the set $\mathcal{N}(k_1, k_2, \dots, k_n)$, therefore it has measure zero.

\end{demo}

\end{document}